\documentclass[12pt]{article}

\usepackage{latexsym,amsfonts,amstext,amsbsy,amsmath,amssymb,amscd,
 graphicx,epsfig,eucal,bbm}
\usepackage{mathrsfs}
\usepackage{pb-diagram}
\usepackage[all]{xy}
\usepackage{makeidx}
\usepackage[romanian]{babel}
\usepackage{indentfirst}
\usepackage{multicol}

\begin{document}
\pagestyle{myheadings}
\begin{center}
{\large \bf{ DYNAMICAL SYSTEMS ON LEIBNIZ ALGEBROIDS}}\\[0.2cm]
\end{center}
\begin{center}
{\bf Gheorghe IVAN and Dumitru OPRI\c S }
\end{center}
\date{}
\def \e{\varepsilon}
\def \w{\widetilde}
\def \b{$\Box$}

{\bf Abstract.} {\small In this paper we study the differential
systems on Leibniz algebroids. We introduce a class of almost
metriplectic manifolds as a special case of Leibniz
 manifolds. Also, the notion of almost metriplectic algebroid is introduced.

 These types of algebroids are used in the presentation of associated differential systems. We give some interesting examples of differential systems on
  algebroids and the orbits of the solutions of corresponding systems are
  established.}
{

\section { Introduction}

Lie algebroids have been introduced repeatedly into differential
geometry since the early $ 1950$'s, and also into physics and
algebra, under a wide variety of names.

It is well known that there exists a one-to-one correspondence
between Lie algebroid structures on a vector bundle $ \pi : E\to M
$ and linear Poisson structures on the dual vector bundle $
\pi^{*} : E^{*}\to M$. This correspondence can be extended to much
wider class of binary operations (brackets) on sections of $ \pi $
on one side, and linear $ 2 $ - contravariant tensor field on $
E^{*} $ on the other side. It is not necessary for these
operations to be antisymmetric or to satisfy the Jacoby identity.
The vector bundle $ \pi $ togheter with a bracket operation, or
the equivalent $ 2 $ - contravariant tensor field, will be called
an algebroid. We mention the concept of Loday algebras, i.e.
Leibniz algebras in the sense of Loday which are " non
antisymmetric Lie algebras ".

Weinstein 's paper on Lagrangian mechanics and groupoids (see [6])
roused new interest into the field of algebroids and groupoids.
Weinstein introduces  " Lagrangian systems " on a Lie algebroid by
means of a Legendre-type map from $ E $ to $ E^{*} $ associated to
a given function $ L $ on $ E $. The local coordinate expression
of such equations reads
$${\dot x}^{i} = \rho_{a}^{i}(x)y^{a}$$
$$\frac{d}{dt}(\frac{\partial L}{\partial y^{a}}) = \rho_{a}^{i}\frac{\partial L}{\partial x^{i}} - C_{ab}^{d}y^{b}\frac{\partial L}{\partial y^{d}} $$
where the $ ( x^{i}) $ are coordinates on $ M $, $ ( y^{a} ) $ are
fibre coordinates on $ E $ and the $ C_{ab}^{d} $ are structure
functions coming from the Lie algebroid structure. Note that, more
generally,  the equations of the form
$${\dot x}^{i} = \rho_{a}^{i}(x) y^{a}$$
$${\dot y}^{a} = f^{a}(x,y)$$
were called " second - order equations on a Lie algebroid " by
Weinstein.

This paper shows that various relevant dynamical systems can be
described as vector fields associated to smooth functions via a
bracket that defines  a Leibniz structure, defined of Leibniz
algebroid.

In this paper we go all the way in this direction and we work with
a bracket, first introduced in (Grabowski and Urbanski, [2]), that
is just required to be linear and a derivation on each its entres.
The derivation property,
 also known as the Leibniz rule, justifies why we refer to this structure as Leibniz bracket.

 This construction should not be mistaken with the Leibniz structures (also called Loday algebras) introduced by Loday ([ 4 ]) in the algebraic context.

The notion of Leibniz algebra is a certain noncommutative version
of a Lie algebra. Leibniz algebroids are bundles with a local
Leibniz algebra structure on its sections and a bundle map to the
tangent bundle (called , anchor), which maps the bracket into the
usual Lie bracket of vector fields.

\section { Leibniz systems. Almost Leibniz systems}

Let $ M $ be a smooth manifold of dimension $ n $ and let $ C^{\infty}(M) $ be the ring of smooth functions on it.
{\it A Leibniz bracket} on $ M $ is a bilinear map $ [\cdot,\cdot ] : C^{\infty}(M)\times C^{\infty}(M) \to C^{\infty}(M) $ such that it is a derivation on each entry,
that is, for all $ f, g, h \in C^{\infty}(M) $ the following relations hold:
\begin{equation}
[ f g , h ] = [ f , h ] g + f [ g , h ] ~~~ \hbox{ and }~~~ [ f , g h ] = g [ f , h ] + [ f, g ] h .\label{1}
\end{equation}

We will say that the pair $ ( M, [ \cdot, \cdot ] ) $ is a {\it Leibniz manifold}. If
the bracket $ [ \cdot, \cdot ] $ is antisymmetric, that is $ [ f , g ] = - [ g , f ], ( \forall ) f,g\in C^{\infty}(M) $, then we say that
$ ( M, [ \cdot, \cdot ] ) $ is an {\it almost Poisson manifold}.

 Let $ ( M, [ \cdot, \cdot ] ) $ be a  Leibniz manifold and let $ h $ be a smooth function of $ M $. From the fact that a Leibniz structure is a derivation there exists the vector field
  $ X_{h} $ on $ M $ uniquely characterized by the relation :
\begin{equation}
X_{h} ( f ) = [ f , h ] ~~~\hbox { for any } ~~~f \in C^{\infty}(M). \label{2}
\end{equation}

We will call $ X_{h} $ the {\it Leibniz vector field} associated to the {\it Hamiltonian function }
$~h \in C^{\infty}(M)$. Note that since  $ [ \cdot, \cdot ] $ is a derivation on each their arguments, they only depend on the first derivatives of the functions and thus,
we can define  the $ 2 $- contravariant tensor field $ B $ on $ M $ given by:
\begin{equation}
B ( d f, d g ) = [ f , g ] , ~~~\hbox{ for all }~~~ f,g\in C^{\infty}(M). \label{3}
\end{equation}

We call {\it dynamical Leibniz system}, the dynamical system associated to the vector field $ X_{h} $.

If $ ( x^{i}), i=\overline{1,n} $ is a system of local coordinates on $ M $, in a local chart
 the dynamical Leibniz system is given by:
\begin{equation}
{\dot x}^{i} = [ x^{i} , h ] , ~~\hbox{ where }~~ [ x^{i} , h ] =
X_{h}^{i} = B ( dx^{i}, d h ) = B^{ij} \frac {\partial h}{\partial
x^{j}}. \label{4}
\end{equation}

Note that the symplectic and Poisson manifolds are particular cases of Leibniz manifolds.
 We now briefly present other non trivial examples.

{\bf Example 1.} $ ( i ) $ {\bf Pseudometric bracket and gradient
dynamical systems}. Let $ {\bf g } $ be a  $ 2 $ - contravariant
non degenerate
 symmetric tensor field on $ M $.
 Let $ [\cdot,\cdot ] : C^{\infty}(M)\times C^{\infty}(M) \to C^{\infty}(M) $  be the Leibniz bracket defined by:
\begin{equation}
[ f , h ] = {\bf g } ( d f, d h ) , ~~\hbox{ for all }~~ f,h \in
C^{\infty}(M). \label{5}
\end{equation}

The bracket given by $ ( 5 ) $ is called the {\it pseudometric bracket associated to } ${\bf g} $. This bracket  is clearly symmetric
and non degenerate and the Leibniz vector field  $ X_{h} $ associated to any function $ h \in C^{\infty}(M) $ is such that
$ X_{h} f = {\bf g } ( d f, d h ) . $  These brackets are also called {\it Beltrami brackets}.

Let on $ {\bf R}^{3} $ the constant $ 2 $ - contravariant tensor
field $ {\bf g} = ( g^{ij} ) $ given by:
\begin{equation}
{\bf g} = \left ( \begin{array}{ccc}
s_{1}\gamma_{1} & 0 & 0 \\
0 & s_{2}\gamma_{2} & 0 \\
0 & 0 & s_{3}\gamma_{3} \\
\end{array}\right )
\label{ 6 }
\end{equation}
where the parameters $ s_{1}, s_{2}, s_{3} \in \{ -1, 1 \} $ and $
\gamma_{1}, \gamma_{2},\gamma_{3} $ are real numbers satisfying
the relation  $ \gamma_{1} + \gamma_{2} + \gamma_{3} = 0 $.  This
system happens to be a particular of pseudometric bracket and the
differential system for
 $ h = x^{1} x^{2} x^{3} $ is given by:
\begin{equation}
{\dot x}^{1} = s_{1}\gamma_{1} x^{2} x^{3},~~~ {\dot x}^{2} = s_{2}\gamma_{2} x^{1} x^{3},~~~ {\dot x}^{3} = s_{3}\gamma_{3} x^{1} x^{2}.
\label{ 7 }
\end{equation}

$ ( ii ) $ {\bf The almost metriplectic systems }. Let  $ P $ be a
$ 2 $ - contravariant antisymmetric tensor field  on $ M $ and $
{\bf g } $ be a  $ 2 $ - contravariant non degenerate symmetric
tensor field on $ M $.
 We define the bracket $ [\cdot,\cdot ] : C^{\infty}(M)\times C^{\infty}(M) \to C^{\infty}(M) $  by:
\begin{equation}
[ f , h ] = P( d f, d h ) + {\bf g } ( d f, d h ) , ~~\hbox{ for all }~~ f,h \in C^{\infty}(M). \label{8}
\end{equation}

This bracket is clearly a Leibniz bracket. We will say that $ (M,
P , {\bf g }, [ \cdot, \cdot ]) $ is an {\it almost metriplectic
manifold}.  If the tensor field $ P $ is Poisson, then  $ (M, P ,
{\bf g }, [ \cdot, \cdot ]) $ is a {\it metriplectic manifold},
see [1].

The differential system associated to the bracket $ (8) $ is given
by:
\begin{equation}
{\dot x}^{i} = [ x^{i} , h ] ,~~\hbox{ where}~~ [ x^{i} , h ] = P
( d x^{i}, d h ) + {\bf g} ( dx^{i}, d h ) = P^{ij} \frac
{\partial h}{\partial x^{j}} + g^{ij} \frac {\partial h}{\partial
x^{j}}. \label{9}
\end{equation}

Let be the $ 2 $ - contravariant tensors fields $ P = (P^{ij} )  ,
{\bf g } = ( g^{ij} ) $ on $ {\bf R}^{3} $ and the function $ h\in
C^{\infty}({\bf R}^{3} ) $ given by:
\begin{equation}
P = \left ( \begin{array}{ccc}
0 & x^{3} & - x^{2}\\
-x^{3} & 0 & x^{1}\\
x^{2} & - x^{1}& 0\\
\end{array}\right )
\label { 10}
\end{equation}
\begin{equation}
{\bf g} =\left ( \begin{array}{ccc}
-a^2_{2}(x^{2})^{2}- a^2_{3}(x^{3})^{2} & a_{1} a_{2} x^{1} x^{2} & a_{1} a_{3} x^{1} x^{3}\\
a_{1} a_{2} x^{1} x^{2} & -a^2_{1}(x^{1})^{2}- a^2_{3}(x^{3})^{2} & a_{2} a_{3} x^{2} x^{3} \\
a_{1} a_{3} x^{1} x^{3} & a_{2} a_{3} x^{2} x^{3} & -a^2_{1}(x^{1})^{2}- a^2_{2}(x^{2})^{2} \\
\end{array}\right )
\label {11}
\end{equation}
\begin{equation}
h = \frac{1}{2} ( a_{1} + 1 )( x^{1})^{2} + \frac{1}{2} ( a_{2} + 1 )( x^{2})^{2} + \frac{1}{2} ( a_{3} + 1 )( x^{3})^{2}.
\label { 12}
\end{equation}

The differential system associated to the bracket $ ( 9 ) $ is
given by:
\begin{equation}
\left\{ \begin{array}{c}
{\dot x}^{1} = ( a_{3} - a_{2} ) x^{2} x^{3} + a_{2}( a_{1} - a_{2} ) x^{1} ( x^{2})^{2} +  a_{3}( a_{1} - a_{3} ) x^{1} ( x^{3})^{2}\\
{\dot x}^{2} = ( a_{1} - a_{3} ) x^{1} x^{3} + a_{3}( a_{2} - a_{3} ) x^{2} ( x^{3})^{2} +  a_{1}( a_{2} - a_{1} ) x^{2} ( x^{1})^{2}\\
{\dot x}^{3} = ( a_{2} - a_{1} ) x^{1} x^{2} + a_{1}( a_{3} - a_{1} ) x^{3} ( x^{1})^{2} +  a_{2}( a_{3} - a_{2} ) x^{3} ( x^{2})^{2}\\
\end{array}\right.
\label { 13 }
\end{equation}

The differential system $ ( 13 ) $ is the {\it revised system } of
the rigid body }, see [ 1 ].\hfill \b

Let $ M $  be a smooth manifold and  $ P $ and ${\bf g } $ be two
$ 2 $ - contravariant tensors fields  on $ M$. We define the
bracket $ [\cdot,( \cdot, \cdot ) ] : C^{\infty}(M)\times
C^{\infty}(M)\times C^{\infty}(M) \longrightarrow C^{\infty}(M) $
by:
\begin{equation}
[ f , ( h_{1}, h_{2} )] = P( d f, d h_{1} ) + {\bf g } ( d f, d h_{2} ) ,~ \hbox{ for all }~ f,h_{1}, h_{2} \in C^{\infty}(M). \label{14}
\end{equation}

{\bf Proposition 1.} {\it The bracket $ ( 14 ) $ satisfy the following relations:
\begin{equation}
[ f_{1} f, ( h_{1}, h_{2} )] = f_{1} [ f, ( h_{1} , h_{2}) ] + f [ f_{1}, ( h_{1} , h_{2}) ] \label{15}
\end{equation}
\begin{equation}
[[ f , l ( h_{1}, h_{2} )] = l [ f , ( h_{1}, h_{2} )] + h_{1} P( d f, d l ) + h_{2} {\bf g } ( d f, d l_{2} ) \label{16}
\end{equation}
for all  $ f, f_{1} , h_{1}, h_{2}, l \in C^{\infty}(M) $.}\hfill
\b

The bracket given by $ ( 14 ) $ is called the {\it almost Leibniz bracket}, and we say that
 $ ( M, P , {\bf g }, [ \cdot, (\cdot,\cdot ) ] ) $ is an {\it almost Leibniz manifold}.

By direct calculation we obtain the following result.

{\bf Proposition 2.} {\it If the tensor field $ P $ is antisymmetric and the tensor field $ {\bf g} $ is symmetric and there exist
 $ h_{1}, h_{2} \in C^{\infty}(M) $ such that $ P(dh_{2},df)=0, ~{\bf g}(dh_{1}, df)=0 $ for all $f\in C^{\infty}(M) $,  then the bracket $( 14 ) $ satisfies the relation:
\begin{equation}
[ f, ( h_{1}, h_{2} )] = [ f, ( h , h ) ] , ~~\hbox{ for}~~ h = h_{1} + h_{2}. \label{17}
\end{equation}
}\hfill \b

 From the fact that the almost Leibniz structure $ (14 ) $ satisfies the relation $ ( 15 ) $ , there exists the vector field $ X_{h_{1} h_{2}} $ on $ M $ uniquelly characterized by the relation :
\begin{equation}
X_{h_{1} h_{2}}(f) = [ f, ( h_{1}, h_{2} )] ~~\hbox{ for any }~~
f\in C^{\infty}(M). \label{18}
\end{equation}

We will call an {\it almost Leibniz dynamical system}, the dynamical system associated of the vector field $ X_{h_{1} h_{2}} $
given by $ ( 18) $.

In  system of local coordinates $ ( x^{i}), i=\overline{1,n} $ on $ M $,
 the almost Leibniz dynamical system is given by:
\begin{equation}
{\dot x}^{i} = [ x^{i} , ( h_{1}, h_{2} ) ] \label{19}
\end{equation}
where
\begin{equation}
[ x^{i} , ( h_{1} , h_{2} ) ] = X_{h_{1}h_{2}}^{i} = P^{ij} \frac {\partial h_{1}}{\partial x^{j}} + g^{ij}\frac{\partial h_{2}}{\partial x^{j}}. \label{20}
\end{equation}

{\bf Example 2.} Let be the $ 2 $ - contravariant tensors fields $
P = (P^{ij} )  , {\bf g } = ( g^{ij} ) $ on $ {\bf R}^{3} $ and
the functions $ h_{1}, h_{2}\in C^{\infty}({\bf R}^{3} ) $ given
by:
\begin{equation}
P = \left ( \begin{array}{ccc}
0 & 1 & 0\\
-1 & 0 & x^{1}\\
0 & - x^{1} & 0\\
\end{array}\right ),~~{\bf g} =\left ( \begin{array}{ccc}
0 & 0 & 0\\
0 & - (x^{3})^{2} & 0 \\
0 &  0 & - (x^{2})^{2} \\
\end{array}\right )\label{21}
\end{equation}
\begin{equation}
h_{1} = \frac{1}{2} ( x^{2})^{2} + \frac{1}{2}( x^{3})^{2}~~,~~ h_{2} = \frac{1}{2}(x^{1})^{2} + x^{3}.
\label { 22}
\end{equation}

The almost Leibniz dynamical system associated to structure $ (
{\bf R}^{3}, P, {\bf g}, [\cdot, (\cdot, \cdot ) ] ) $ for the
above tensors fields and functions is given by:
\begin{equation}
{\dot x}^{1} = x^{2},~~ {\dot x}^{2} = x^{1} x^{3}, ~~ {\dot x}^{3} = - x^{1} x^{2} - (x^{2})^{2}.
\label{23}
\end{equation}
 \hfill \b

{\bf Example 3. } If we take  the tensors fields $ P = (P^{ij} ) $
and  $ {\bf g } = ( g^{ij} ) $ on $ {\bf R}^{3} $ and the
functions $ h_{1}, h_{2}\in C^{\infty}({\bf R}^{3} ) $ given by:
\begin{equation}
P = \left ( \begin{array}{ccc}
0 & -x^{3} & x^{2}\\
x^{3} & 0 & 0\\
- x^{2} & 0 & 0\\
\end{array}\right ) ~~,~~{\bf g} =\left ( \begin{array}{ccc}
- x^{3} & 0 & 0\\
0 & 0 & 0 \\
0 &  0 & - x^{1} \\
\end{array}\right ) \label{24}
\end{equation}
\begin{equation}
h_{1} = \frac{1}{2} ( x^{1})^{2} + x^{3}~,~~ h_{2} = \frac{1}{2}(x^{2})^{2} +\frac{1}{2}( x^{3})^{2}
\label { 25}
\end{equation}
then the almost Leibniz dynamical system associated to
 $( {\bf R}^{3}, P, {\bf g}, [\cdot, (\cdot, \cdot ) ] ) $ for
the above tensors fields and functions is given by:
\begin{equation}
{\dot x}^{1} = x^{2},~~ {\dot x}^{2} = x^{1} x^{3}, ~~ {\dot x}^{3} = - x^{1} x^{2} - x^{1} x^{3}.
\label{26}
\end{equation}
\hfill \b

\section { Leibniz algebroids}

Let $ M $ be a smooth manifold of dimension $ n $, let $ \pi : E
\to M $ be a vector bundle and $  \pi^{*} : E^{*} \to M $ the dual
vector bundle. By $~\Sigma = \Gamma(M,E) $ we denote the sections
of $ \pi $.

A {\it Leibniz algebroid structure }({\it pseudo - Lie algebroid
structure }) on a vector bundle $~ \pi : E \to M $ is
 given by a bracket ( bilinear operation ) $ [ \cdot, \cdot ] $ on the space of sections $ \Sigma $ and two vector bundle morphisms
$ \rho_{1}, \rho_{2} : E \to TM $ (called the {\it left} and the
{\it right anchor}, respectively) such that
\begin{equation}
[ f \sigma_{1}, g \sigma_{2} ] = f \rho_{1}(\sigma_{1})(g)\sigma_{2} - g \rho_{2}(\sigma_{2})(f) \sigma_{1} + f g [ \sigma_{1}, \sigma_{2}] \label { 27}
\end{equation}
for all $ \sigma_{1}, \sigma_{2}\in \Sigma $ and $ f, g \in C^{\infty}(M).$

 A vector bundle  $ \pi : E \to M $ endowed with a Leibniz algebroid structure $ ([\cdot, \cdot ], \rho_{1}, \rho_{2}) $ on $ E $ , is called {\it
 Leibniz algebroid} over $ M$ and denoted by $(E, [\cdot, \cdot ], \rho_{1}, \rho_{2}).$

 A Leibniz algebroid with an antisymmetric bracket $ [\cdot, \cdot ] $ ( in this case we have $~\rho_{1}= \rho_{2} $ ) is called {\it pre - Lie algebroid}.
 A Leibniz algebroid with a symmetric bracket $ [\cdot, \cdot ] $ ( in this case we have $~\rho_{2}= -\rho_{1} $ ) is called {\it symmetric algebroid}.

In the following, we establish a correspondence between Leibniz algebroid structures on the vector bundle $ \pi : E \to M $ and the $ 2 $- contravariant tensor fields on
bundle manifold $ E^{*} $ of the dual vector bundle $ \pi^{*} : E^{*} \to M $.

For a given section $ \sigma \in \Sigma ,$ we define a function $
i_{E^{*}}\sigma $ on $ E^{*} $ by the relation:
\begin{equation}
i_{E^{*}}\sigma(a) = < \sigma(\pi^{*}(a)),a > ,~~\hbox { for } ~~  a\in E^{*},  \label {28}
\end{equation}
where $ <\cdot, \cdot > $ is the canonical pairing between $ E $ and $ E^{*}. $

Let $ \Lambda $ be a $ 2 $ - contravariant tensor field on $ E^{*} $ and the bracket $ [\cdot, \cdot ]_{\Lambda} $ of functions given by the relation :
 \begin{equation}
 [ f, g ]_{\Lambda} = \Lambda( df, d g )\label {29}
 \end{equation}
for all $f,g\in C^{\infty}(E^*).$

 For  a given $ 2 $ - contravariant tensor field   $ \Lambda $ on $ E^{*}, $ we say that $ \Lambda $
 is {\it linear}, if for each pair $ ( \mu_{1}, \mu_{2} ) $ of sections of $ \pi^{*} $ , the function
 $ \Lambda( d i_{E^{*}}\mu_{1}, d i_{E^{*}}\mu_{2} ) $ defined on $ E^{*} $ is linear.

{\bf Theorem 1.} (Grabowski and Urbanski, 1997, [2]) {\it  For
every Leibniz algebroid structure on $ \pi : E \to M $
 with the bracket $ [ \cdot, \cdot ] $ and the anchors $ \rho_{1}, \rho_{2} $ , there exists an unique $ 2 $ - contravariant tensor field
$ \Lambda $ on $ E^{*} $ such that the following relations hold:
\begin{equation}
i_{E^{*}}[ \sigma_{1}, \sigma_{2} ] = [ i_{E^{*}}\sigma_{1}, i_{E^{*}}\sigma_{2}  ]_{\Lambda}\label { 30}
\end{equation}
\begin{equation}
\pi^{*}(\rho_{1}(\sigma)(f))= [ i_{E^{*}}\sigma ,\pi^{*} f]_{\Lambda}~,~ \pi^{*}(\rho_{2}(\sigma)(f))= [ \pi^{*} f, i_{E^{*}}\sigma ]_{\Lambda}  \label { 31}
\end{equation}

for all $ \sigma, \sigma_{1}, \sigma_{2} \in \Sigma $ and $ f\in C^{\infty}(M). $

Conversely, every $ 2 $ - contravariant linear tensor field $
\Lambda $ on $ E^{*} $ defines a Leibniz algebroid on $ E $ by the
relations $ ( 30 ) $ and $ ( 31 ) $.}\hfill \b

Let $ ( x^{i} ), i=\overline{1,n} $ be a local coordinate system
on $ M $ and let $ \{ e_{1}, \ldots, e_{m} \} $ be a basis of
local sections of $ E $ ( $ dim ~M = n, dim~ E = m $ ). We denote
by $ \{ e^{1}, \ldots, e^{m}\} $ the dual basis of local sections
of $ E^{*} $ and $ ( x^{i}, y^{a} ) $ ( resp., $ ( x^{i}, \xi_{a}
) $ ) the corresponding coordinate system on $ E $ ( resp., $
E^{*} $ ).

It is easy to see that every linear $ 2 $ - contravariant tensor
field $ \Lambda $ on $ E^{*} $ has the form:
\begin{equation}
\Lambda = C_{ab}^{d}\xi_{d}\frac{\partial }{\partial \xi_{a}}\otimes \frac{\partial }{\partial \xi_{b}}   +
\rho_{1 a}^{i}\frac{\partial }{\partial \xi_{a}}\otimes \frac{\partial }{\partial x^{i}} -
\rho_{2 a}^{i}\frac{\partial }{\partial x^{i}}\otimes \frac{\partial }{\partial \xi_{a}} , \label { 32}
\end{equation}
where $ C_{ab}^{d}, \rho_{1 a}^{i}, \rho_{2 a}^{i}\in
C^{\infty}(M) $  are functions of $ (x^{i}). $

The correspondence  between $ \Lambda $ and a Leibniz algebroid
structure is given by the following relations:
\begin{equation}
[ e_{a}, e_{b}] = C_{ab}^{d}e_{d}~,~ \rho_{1}(e_{a}) = \rho_{1 a}^{i}\frac{\partial}{\partial x^{i}}~,~
\rho_{2}(e_{a}) = \rho_{2 a}^{i}\frac{\partial}{\partial x^{i}}.\label{33}
\end{equation}

We call a {\it dynamical system on Leibniz algebroid $ \pi : E\to M $ }, the dynamical system associated to vector field $ X_{h} $ with
$ h\in C^{\infty}(M) $ given by:
\begin{equation}
X_{h}(f) = \Lambda(df, dh),~\hbox { for all }~ f\in C^{\infty}(M).\label { 34}
\end{equation}

In a system of local coordinates, the dynamical system $ ( 34) $
is given by:
\begin{equation}
{\dot \xi}_{a} = [ \xi_{a}, h]_{\Lambda},~~~{\dot x}^{i} = [ x^{i}, h ]_{\Lambda}\label {35}
\end{equation}
where
\begin{equation}
[ \xi_{a}, h]_{\Lambda}= C_{ab}^{d}\xi_{d}\frac{\partial h}{\partial \xi_{b}} + \rho_{1 a}^{i}\frac{\partial h}{\partial x^{i}}~~,~~
[ x^{i}, h ]_{\Lambda} = - \rho_{2 a}^{i}\frac{\partial h}{\partial \xi_{a}} .\label{36}
\end{equation}

If the Leibniz algebroid is a pre - Lie algebroid, then the
dynamical system $ (35 )$ and $ ( 36) $ is given by:
\begin{equation}
{\dot \xi}_{a} = C_{ab}^{d}\xi_{d}\frac{\partial h}{\partial \xi_{b}} + \rho_{1 a}^{i}\frac{\partial h}{\partial x^{i}} ~~,~~
{\dot x}^{i}  = - \rho_{1 a}^{i}\frac{\partial h}{\partial \xi_{a}}, ~ C_{ab}^{d} = -C_{ba}^{d}.\label{37}
\end{equation}

If the Leibniz algebroid is a symmetric algebroid, then the
dynamical system $ (35)$ and $ ( 36 ) $ is given by:
\begin{equation}
{\dot \xi}_{a} = C_{ab}^{d}\xi_{d}\frac{\partial h}{\partial \xi_{b}} + \rho_{1 a}^{i}\frac{\partial h}{\partial x^{i}}~~,~~
{\dot x}^{i}  = \rho_{1 a}^{i}\frac{\partial h}{\partial \xi_{a}}, ~ C_{ab}^{d} =  C_{ba}^{d}.\label{38}
\end{equation}

{\bf Example 4.} Let the vector bundle $ \pi : E = {\bf R}^{3}\times {\bf R}^{3}\to {\bf R}^{3} $ and
 $ \pi^{*} : E^{*}={\bf R}^{3}\times ({\bf R}^{3})^{*}\to {\bf R}^{3} $ the dual vector bundle. We consider on $ E^{*} $ the $ 2 $ - contravariant linear
 tensor field $ \Lambda $, the anchors $ \rho_{1}, \rho_{2} $ and the function $ h $ given by:
\begin{equation}
P = \left ( \begin{array}{ccc}
0 & -\xi_{3}x^{3} & \xi_{2}x^{2}\\
\xi_{3}x^{3} & 0 & -\xi_{1}x^{1}\\
- \xi_{2}x^{2} & \xi_{1}x^{1} & 0\\
\end{array}\right ),~~\rho_{1} =\left ( \begin{array}{ccc}
0 &- x^{3} & x^{2}\\
x^{3} & 0 & 0 \\
- x^{2} & 0 & 0 \\
\end{array}\right )
\label{39}
\end{equation}
\begin{equation}
\rho_{2} =\left ( \begin{array}{ccc}
0 & -1 & 0\\
1 & 0 & -x^{1} \\
0 & x^{1} & 0 \\
\end{array}\right ),~~\hbox{ and }~~h(x,\xi) =  x^{2}\xi_{2} + x^{3}\xi_{3}.
\label{40}
\end{equation}
The dynamical system  $ ( 35 ) $ and $ ( 36 ) $ associated to Leibniz algebroid $ ( {\bf R}^{3}\times {\bf R}^{3}, P, {\bf g}, \rho_{1}, \rho_{2} ) $ for
 $ h = x^{2}\xi_{2}+ x^{3}\xi_{3} $ is given by:
\begin{equation}
\left\{ \begin{array}{ccc} {\dot \xi}_{1} & = &
x^{3}(x^{2}-1)\xi_{2}- x^{2}x^{3}\xi_{3}\cr {\dot \xi}_{2} & = & -
x^{3}x^{1}\xi_{1}\cr {\dot \xi}_{3} & = & x^{1}x^{2}\xi_{1}\cr
\end{array}\right.,~~ \left\{\begin{array}{ccc}
{\dot x}^{1} & = & x^{2}\cr {\dot x}^{2} & = & x^{1}x^{3}\cr {\dot
x}^{3} & = & - x^{1}x^{2}\cr
\end{array}\right.
\label{41}
\end{equation}

The dynamical system $ ( 41 )$ is the {\it dynamical system
associated to the Maxwel-Bloch system}. The orbits of the system $
( 41 ) $ represented in the coordinate systems $ O x^{1} x^{2}
x^{3} $ and $ O \xi_{1} \xi_{2} \xi_{3} $ are given in the figures
Fig.1 and Fig. 2.

\begin{center}\begin{tabular}{cc}
\epsfxsize=6cm \epsfysize=5cm
 \epsffile{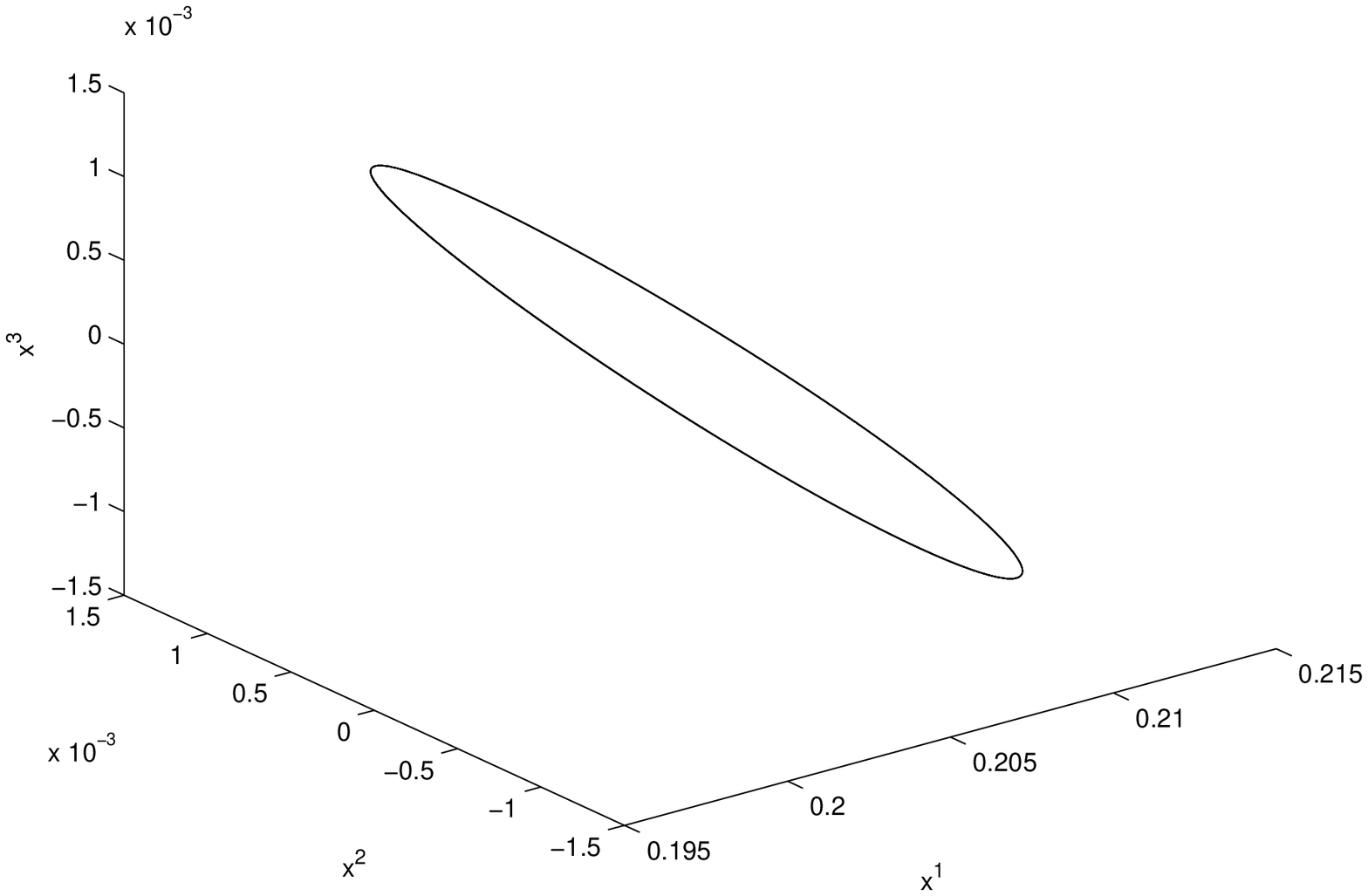} &
\epsfxsize=6cm \epsfysize=5cm \epsffile{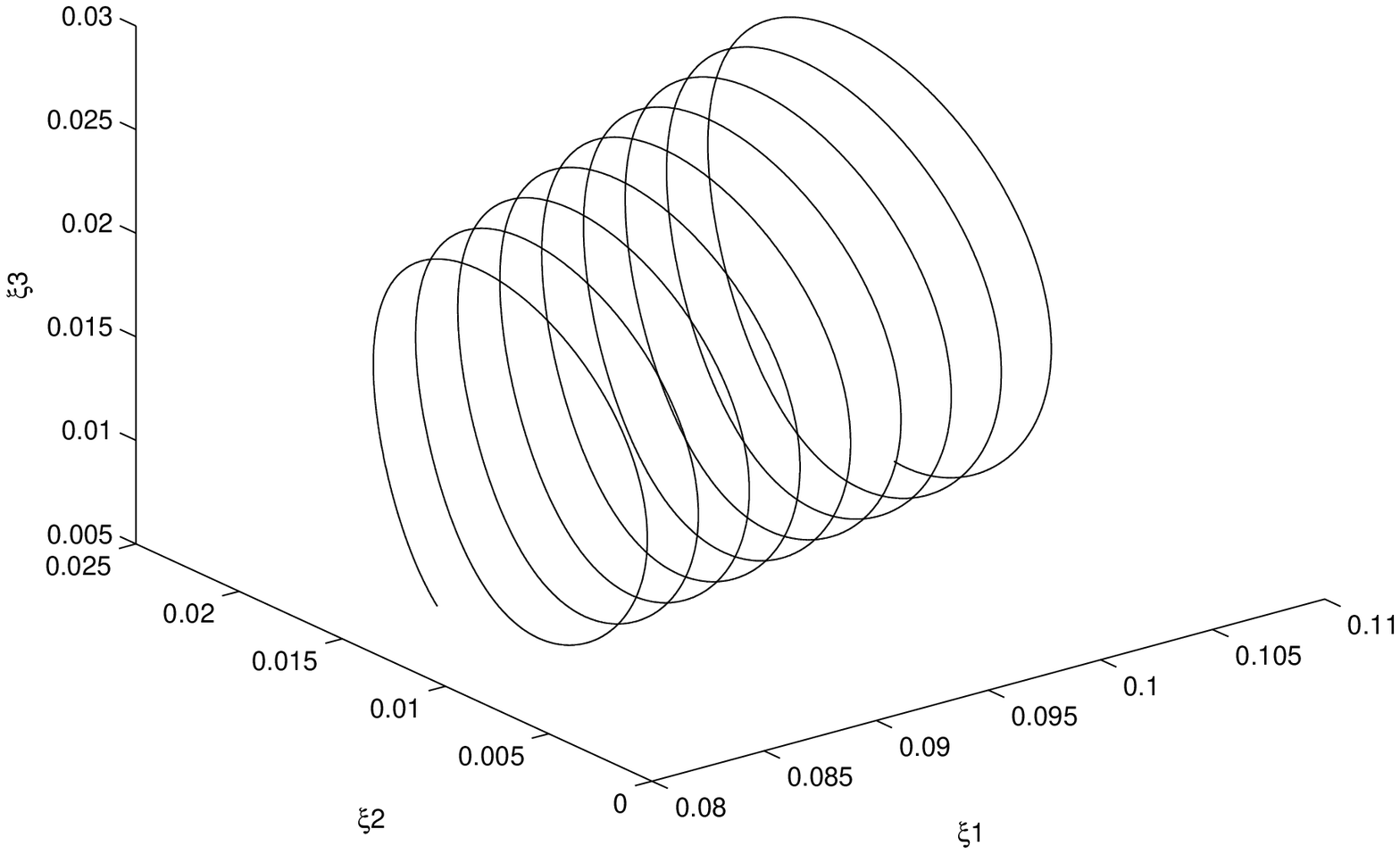}\\
Fig. 1 & Fig. 2
\end{tabular}
\end{center}

\section { Almost metriplectic algebroids}

Let $ \pi : E\to M$ be a vector bundle , an antisymmetric bracket $ [ \cdot, \cdot ]_{1}$ on $ E $ and
$ [\cdot, \cdot ]_{2} $ a symmetric bracket on $ E.$ A Leibniz algebroid endowed with the brackets
$ [\cdot, \cdot ]_{1} $ and $ [\cdot, \cdot ]_{2} $ is called {\it almost metriplectic algebroid.}

Let $ \Lambda_{1} $ and $  \Lambda_{2} $ the $ 2 $ - contravariant tensors fields on the bundle manifold $ E ^{*} $ associated with
brackets  $ [\cdot, \cdot ]_{1} $ and $ [\cdot, \cdot ]_{2}. $

We consider the bracket $ [\cdot , (\cdot,\cdot ) ] : C^{\infty}(E^{*})\times (C^{\infty}(E^{*})\times C^{\infty}(E^{*}))\to C^{\infty}(E^{*}) $ given by:
\begin{equation}
[ f, ( h_{1}, h_{2} )] = \Lambda_{1} ( d f, d h_{1}) + \Lambda_{2}
( d f, d h_{2}) , ~\hbox{ for all }~  f, h_{1}, h_{2} \in
C^{\infty}(E^{*}) . \label{42}
\end{equation}

It is not hard to prove the following proposition.

{\bf Proposition 3.} {\it The bracket given by $ (42) $ satisfy
the relations:
\begin{equation}
[ f_{1}f , ( h_{1}, h_{2} )] =  f_{1} [ f , ( h_{1}, h_{2} )] + f [ f_{1} , ( h_{1}, h_{2} )]\label{43}
\end{equation}
\begin{equation}
[ f , l ( h_{1}, h_{2} )] =  l [ f , ( h_{1}, h_{2} )] + h_{1}\Lambda_{1}(df, dl) +  h_{2}\Lambda_{2}(df, dl) \label{44}
\end{equation}
for all  $ f,f_{1}, h_{1}, h_{2}, l \in C^{\infty}(E^{*})
$.}\hfill \b

From $ (43 )$ it follows that $ X_{h_{1}h_{2}}: C^{\infty}(E^{*})\times C^{\infty}(E^{*})\to C^{\infty}(E^{*})$ given by:
\begin{equation}
X_{ h_{1} h_{2}}(f) =  [ f , ( h_{1}, h_{2} )] ,~~\hbox { for all
} ~~ f\in C^{\infty}(E^{*})  \label{45}
\end{equation}
is a vector field on $ E^{*} $ associated to the functions $ h_{1}, h_{2} \in C^{\infty}(E^{*}). $

We will call an { \it almost metriplectic dynamical system}, the dynamical system associated to the vector field $ X_{h_{1} h_{2}} $
 given by $ ( 45 ) $ with the bracket defined by $ ( 42 ) $.

In a system of local coordinates on $ E^{*} $,  the almost
metriplectic dynamical system $ ( 45) $ is given by:
\begin{equation}
{\dot \xi}_{a} = X_{h_{1} h_{2}}(\xi_{a}) ~,~~{\dot x}^{i} = X_{h_{1} h_{2}}( x^{i} ).\label {46}
\end{equation}
where
\begin{equation}
X_{h_{1} h_{2}}(\xi_{a})=  \Lambda_{1}( d\xi_{a}, dh_{1}) + \Lambda_{2}( d\xi_{a}, dh_{2})~, ~~
X_{h_{1} h_{2}}(x^{i})=  \Lambda_{1}( dx^{i}, dh_{1}) + \Lambda_{2}( dx^{i}, dh_{2})
\label {47}
\end{equation}

Then the almost metriplectic dynamical system $ ( 46) $ and $ ( 47
) $ is given by:
\begin{equation}
\left\{ \begin{array}{c}
{\dot \xi}_{a} = C_{1 ab}^{d}\xi_{d}\frac{\partial h_{1}}{\partial \xi_{b}}+ C_{2 ab}^{d}\xi_{d}\frac{\partial h_{2}}{\partial \xi_{b}} +
\rho_{1 a}^{i}\frac{\partial h_{1}}{\partial x^{i}} + \rho_{2 a}^{i}\frac{\partial h_{2}}{\partial x^{i}}\\
{\dot x}^{i} = -\rho_{1 a}^{i}\frac{\partial h_{1}}{\partial \xi_{a}} + \rho_{2 a}^{i}\frac{\partial h_{2}}{\partial \xi_{a}}.\\
\end{array}\right.
 \label {48}
\end{equation}

{\bf Proposition 4.} {\it Let $ \Lambda_{1} $ and $ \Lambda_{2} $ be two $ 2 $ - contravariant tensors fields on
 $ E^{*} $ for an almost metriplectic algebroid.

$ ( i ) $ If there exist $ h_{1}, h_{2} \in C^{\infty}(E^{*}) $ such that
 $ \Lambda_{1}(d h_{2},df ) = 0 $ and $ \Lambda_{2}(d h_{1},df ) = 0 $ for all $ f\in C^{\infty}(E^{*}) $, then for $ h = h_{1}+ h_{2}~ $ we have
$~[ f , ( h, h )] =  [ f , ( h_{1}, h_{2} )]. $

$(ii) $ If the function $ h_{1} $ is linear on the fibre of the vector bundle $ \pi^{*} : E^{*}\to M $, then there exists $ \Lambda_{2} $ with property that it is determined by $ h_{1}.$}

{\bf Proof.}  $ ( i ) $ From  hypothesis $ \Lambda_{1}(d h_{2},df ) = 0 $ and $ \Lambda_{2}(d h_{1},df ) = 0 $ for all $ f\in C^{\infty}(E^{*}) $ and from $ ( 42 ) $ it follows that\\
$ [f, (h,h)]=\Lambda_{1}(df,dh)+\Lambda_{2}(df,dh)= \Lambda_{1}(df,dh_{1})+ \Lambda_{2}(df,dh_{2}) = [f, ( h_{1}, h_{2})]. $

$ ( ii ) $ From  $\Lambda_{2}(df,dh_{1})= 0 $ for all $ f\in C^{\infty}(E^{*}) $ , we obtain:
\begin{equation}
C_{2 ab}^{d}\xi_{d}\frac{\partial h_{1}}{\partial \xi_{b}}+ \rho_{2 a}^{i}\frac{\partial h_{1}}{\partial x^{i}} = 0 ~~,~~ \rho_{2 a}^{i}\frac{\partial h_{1}}{\partial \xi_{a}} = 0. \label {49}
\end{equation}

If $~ h_{1}(x,\xi)= \xi_{a}h_{1}^{a}(x)$ from $ (49 ) $ follows
\begin{equation}
C_{2 ab}^{d} h_{1}^{b} + \rho_{2 a}^{i}\frac{\partial h_{1}^{d}}{\partial x^{i}} = 0 ~~,~~ \rho_{2 a}^{i} h_{1}^{a}= 0. \label {50}
\end{equation}

Then, it is easy to verify that
$$\left\{\begin{array}{c}
\rho_{2 a}^{i}=\frac{\partial h_{1}^{a}}{\partial x^{i}} ~\hbox{ if } ~ a\neq i ,~~\rho_{2 i}^{a} = - \frac{1}{h_{1}^{i}}
\sum_{ k=1\atop k\neq i}^{n} (\frac{\partial h_{1}^{a}}{\partial x^{k}})^{2},~\hbox { if } a=i \\
C_{2 ab}^{d}= 0 ~\hbox{ if } ~ a\neq b ,~~C_{2 ab}^{d} = - \frac{1}{h_{1}^{b}}\rho_{2 a}^{i} \frac{\partial h^{d}}{\partial x^{i}},~\hbox { if } a=b \\
\end{array}\right.$$
satisfy the relations $ ( 50 ). $\hfill \b

{\bf Example 5.} {\it The dynamical system of the rigid body on the Leibniz algebroid } $ \pi : {\bf R}^{3}\times {\bf R}^{3}\to {\bf R}^{3}.$

The dynamical system of the rigid body on ${\bf R}^{3}$ is given
by:

$${\dot x}^{i} = p^{ij} \frac{\partial h}{\partial x^{j}} , ~~i, j = 1,2,3 $$
where $ ~~h = \frac{1}{2}( a_{1}(x^{1})^{2}+a_{2}(x^{2})^{2}+a_{3}(x^{3})^{2}),~ a_{1} > a_{2} > a_{3} > 0~ $ and \\
$$p = ( p^{ij}) =\left ( \begin{array}{ccc}
0 & x^{3} & -x^{2}\\
-x^{3} & 0 & x^{1}\\
x^{2} & -x^{1} & 0\\
\end{array}\right ) . $$

Let $ \pi^{*} : E^{*}={\bf R}^{3}\times ({\bf R}^{3})^{*}\to {\bf R}^{3}$ the dual vector bundle of $ \pi : E \to {\bf R}^{3} $ and
$ ( x^{i}, \xi_{i} ) $ be a coordinate system on $ E^{*}$. Let the $ 2 $- contravariant tensor field $ \Lambda_{1} $ on $ E^{*} $ given by
 $ ( C_{1},-\rho_{1}, \rho_{1}) $ where $ \rho_{1}= p $ and
$$ C_{1} =\left ( \begin{array}{ccc}
0 & \xi_{3}x^{3} & -\xi_{2}x^{2}\\
-\xi_{3}x^{3} & 0 & \xi_{1}x^{1}\\
\xi_{2}x^{2} & -\xi_{1}x^{1} & 0\\
\end{array}\right ).$$

 The dynamical system for $ \Lambda_{1} $ and $ h_{1}(x,\xi) = a_{1}x^{1}\xi_{1}+ a_{2}x^{2}\xi_{2}+a_{3}x^{3}\xi_{3}$ is
  given by:
\begin{equation}
\left\{ \begin{array}{ccc}
{\dot \xi}_{1} & = & a_{2}\xi_{3}x^{2}x^{3}-a_{3}\xi_{2}x^{2}x^{3}-a_{2}\xi_{2} x^{3}+ a_{3}\xi_{3} x^{2}\cr
{\dot \xi}_{2} & = & a_{3}\xi_{1}x^{1}x^{3}-a_{1}\xi_{3}x^{1}x^{3}-a_{3}\xi_{3} x^{1}+ a_{1}\xi_{1} x^{3} \cr
{\dot \xi}_{3} & = & a_{1}\xi_{2}x^{1}x^{2}-a_{2}\xi_{1}x^{1}x^{2}-a_{1}\xi_{1} x^{2}+ a_{2}\xi_{2} x^{1}\cr
{\dot x}^{1} & = & ( a_{3} - a_{2})x^{2}x^{3}\cr
{\dot x}^{2} & = & ( a_{1} - a_{3})x^{1}x^{3}\cr
{\dot x}^{3} & = & ( a_{2} - a_{1})x^{1}x^{2}\cr
\end{array}\right.
\label{51}
\end{equation}

 The orbits of the dynamical system $ ( 51) $ represented in the coordinate systems $ O x^{1} x^{2} x^{3} $ and $ O \xi_{1} \xi_{2} \xi_{3} $ are given in the figures
Fig. 3 and Fig. 4.

\begin{center}\begin{tabular}{cc}
\epsfxsize=6cm \epsfysize=5cm
 \epsffile{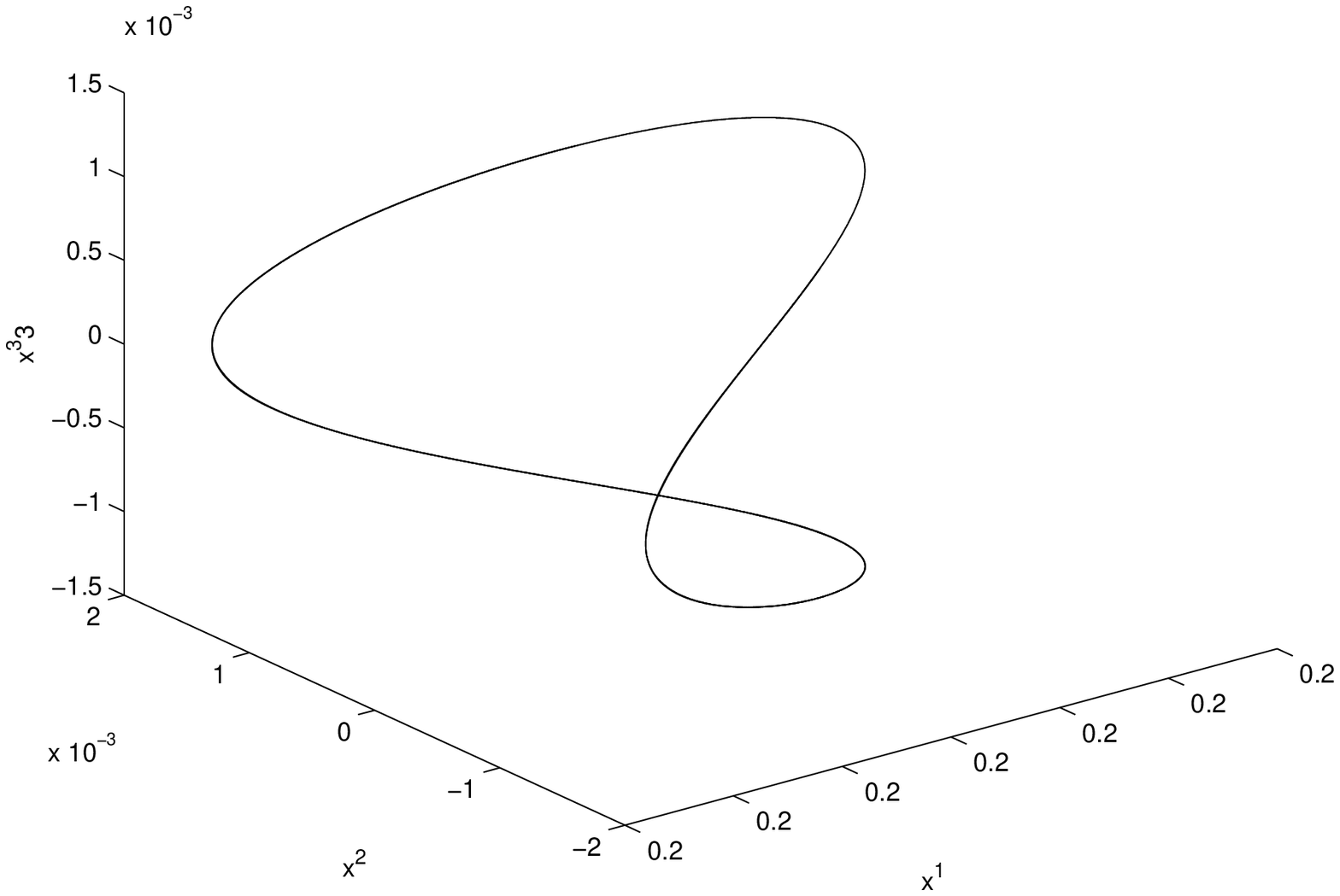} &
\epsfxsize=6cm \epsfysize=5cm \epsffile{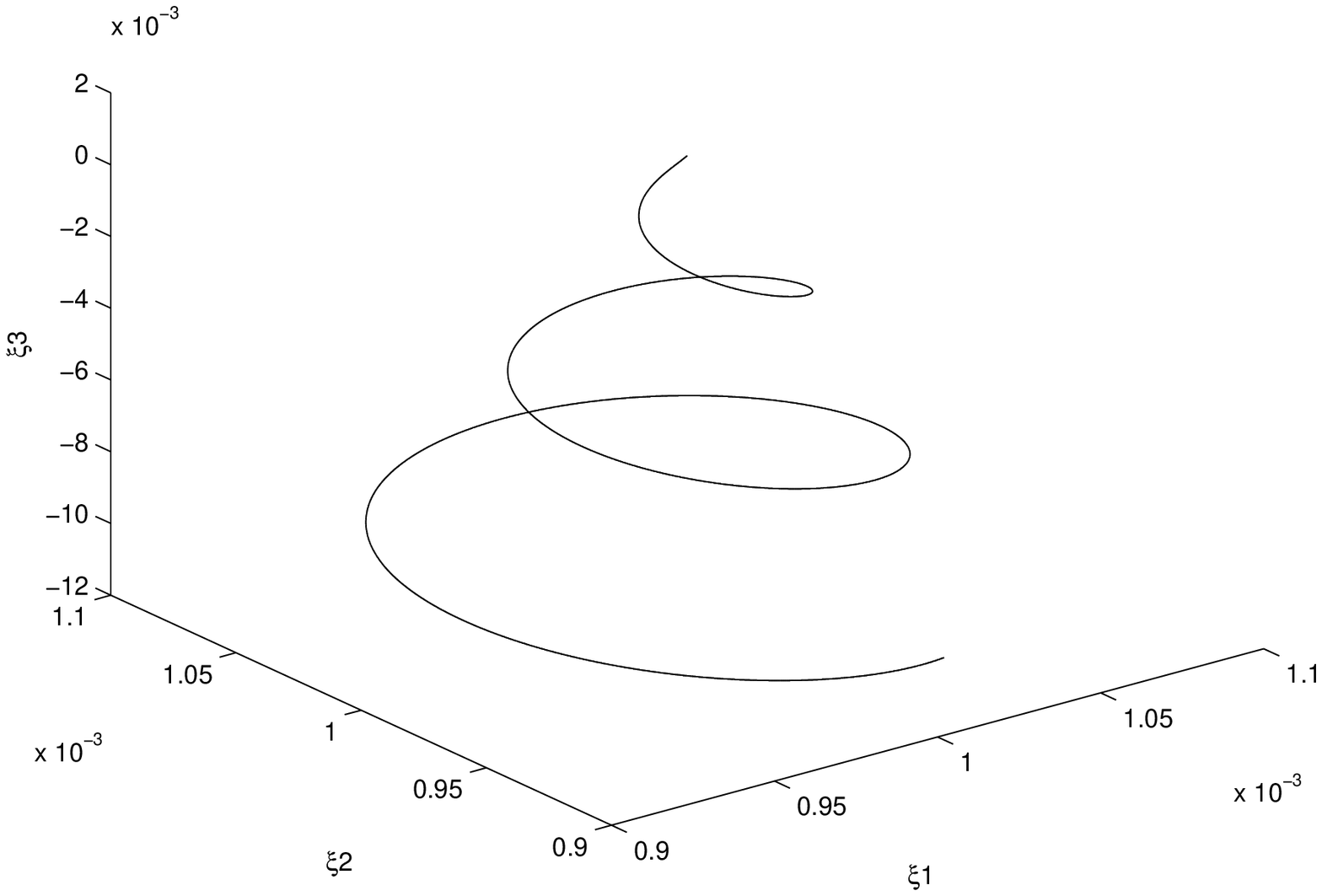}\\
Fig. 3 & Fig. 4
\end{tabular}
\end{center}

The tensor field $ \Lambda_{1} $ satisfies the relation $
\Lambda_{1}(d h_{2}, df ) = 0 $ for all $ f\in C^{\infty}(E^{*})$
and
 $ h_{2}(x,\xi ) = x^{1}\xi_{1}+ x^{2}\xi_{2}+x^{3}\xi_{3}. $

 Applying the Proposition 4 , it follows that there exists
 the $ 2 $ - contravariant tensor field $ \Lambda_{2} $ with property $ \Lambda_{2}(dh_{1}, df ) = 0, $ for all $ f\in C^{\infty}(E^{*})$.

 The $ 2 $ - contravariant tensor field $ \Lambda_{2} $ on $ E^{*} $  is given
 by:
 $ ( C_{2},\rho_{2}, \rho_{2}) $, where
$$C_{2} =\left ( \begin{array}{ccc}
\frac{V_{1}}{x^{1}}\xi_{1}- W_{1} & 0 &0\\
 0 & \frac{V_{2}}{x^{2}}\xi_{2}- W_{2} & 0\\
0 & 0 & \frac{V_{3}}{x^{3}}\xi_{3}- W_{3} \\
\end{array}\right )~\hbox{ and }~~$$
$$\rho_{2} =\left ( \begin{array}{ccc}
- V_{1} & a_{1}a_{2}x^{1}x^{2} &  a_{1} a_{3}x^{1}x^{3} \\
a_{1}a_{2}x^{1}x^{2} & - V_{2} & a_{2} a_{3}x^{2}x^{3} \\
a_{1} a_{3}x^{1}x^{3} & a_{2} a_{3}x^{2}x^{3}  & - V_{3} \\
\end{array}\right ) $$
with\\
$V_{1} = (a_{2})^{2}(x^{2})^{2} +(a_{3})^{2}(x^{3})^{2}, ~~~ W_{1} = (a_{2})^{2}x^{2}\xi_{2}+(a_{3})^{2}x^{3}\xi_{3} , $\\
$V_{2} = (a_{1})^{2}(x^{1})^{2} +(a_{3})^{2}(x^{3})^{2},~~~ W_{2} = (a_{1})^{2}x^{1}\xi_{1} +(a_{3})^{2}x^{3}\xi_{3} , $\\
$V_{3} = (a_{1})^{2}(x^{1})^{2}+ (a_{2})^{2}(x^{2})^{2}, ~~~ W_{3}= (a_{1})^{2}x^{1}\xi_{1} +(a_{2})^{2}x^{2}\xi_{2} .$

 The dynamical system  $ ( 46 ) $ and $ ( 47 ) $ is the following:
\begin{equation}
\left \{ \begin{array}{c}
{\dot \xi}_{1} = (a_{2}( a_{1} - a_{2})x^{1}x^{2}- a_{3}x^{2}x^{3}-a_{2}x^{3})\xi_{2} + (a_{3}( a_{1} - a_{3})x^{1}x^{3}+ a_{2}x^{2}x^{3}+ a_{3}x^{2})\xi_{3}\\
{\dot \xi}_{2} = (a_{1}( a_{2} - a_{1})x^{1}x^{2}+ a_{3}x^{1}x^{3}+a_{1}x^{3})\xi_{1} + (a_{3}( a_{2} - a_{3})x^{2}x^{3}- a_{1}x^{1}x^{3}- a_{3}x^{1})\xi_{3} \\
{\dot \xi}_{3} = (a_{2}( a_{3} - a_{2})x^{2}x^{3} + a_{1}x^{1}x^{2} +a_{2}x^{1})\xi_{2} + (a_{1}( a_{3} - a_{1})x^{1}x^{3}- a_{2}x^{1}x^{2}- a_{1}x^{2})\xi_{1} \\
{\dot x}^{1} = ( a_{3} - a_{2})x^{2}x^{3} + a_{2}( a_{1} - a_{2} )x^{1}(x^{2})^{2} + a_{3}( a_{1}- a_{3} ) x^{1} (x^{3})^{2}\\
{\dot x}^{2} = ( a_{1} - a_{3})x^{1}x^{3} + a_{3}( a_{2} - a_{3} )x^{2}(x^{3})^{2} + a_{1}( a_{2}- a_{1} ) x^{2} (x^{1})^{2} \\
{\dot x}^{3} = ( a_{2} - a_{1})x^{1}x^{2} + a_{1}( a_{3} - a_{1} )x^{3}(x^{1})^{2} + a_{2}( a_{3}- a_{2} ) x^{3} (x^{2})^{2} \\
\end{array}\right.
\label{52}
\end{equation}

 The orbits of the dynamical system $ ( 52) $ represented in the coordinate systems $ O x^{1} x^{2} x^{3} $ and $ O \xi_{1} \xi_{2} \xi_{3} $
 for $ a_{1}=0.6, a_{2} =0.4 , a_{3} =0.2 $  are given in the figures
Fig.5 and Fig. 6.

\begin{center}\begin{tabular}{cc}
\epsfxsize=6cm \epsfysize=5cm
 \epsffile{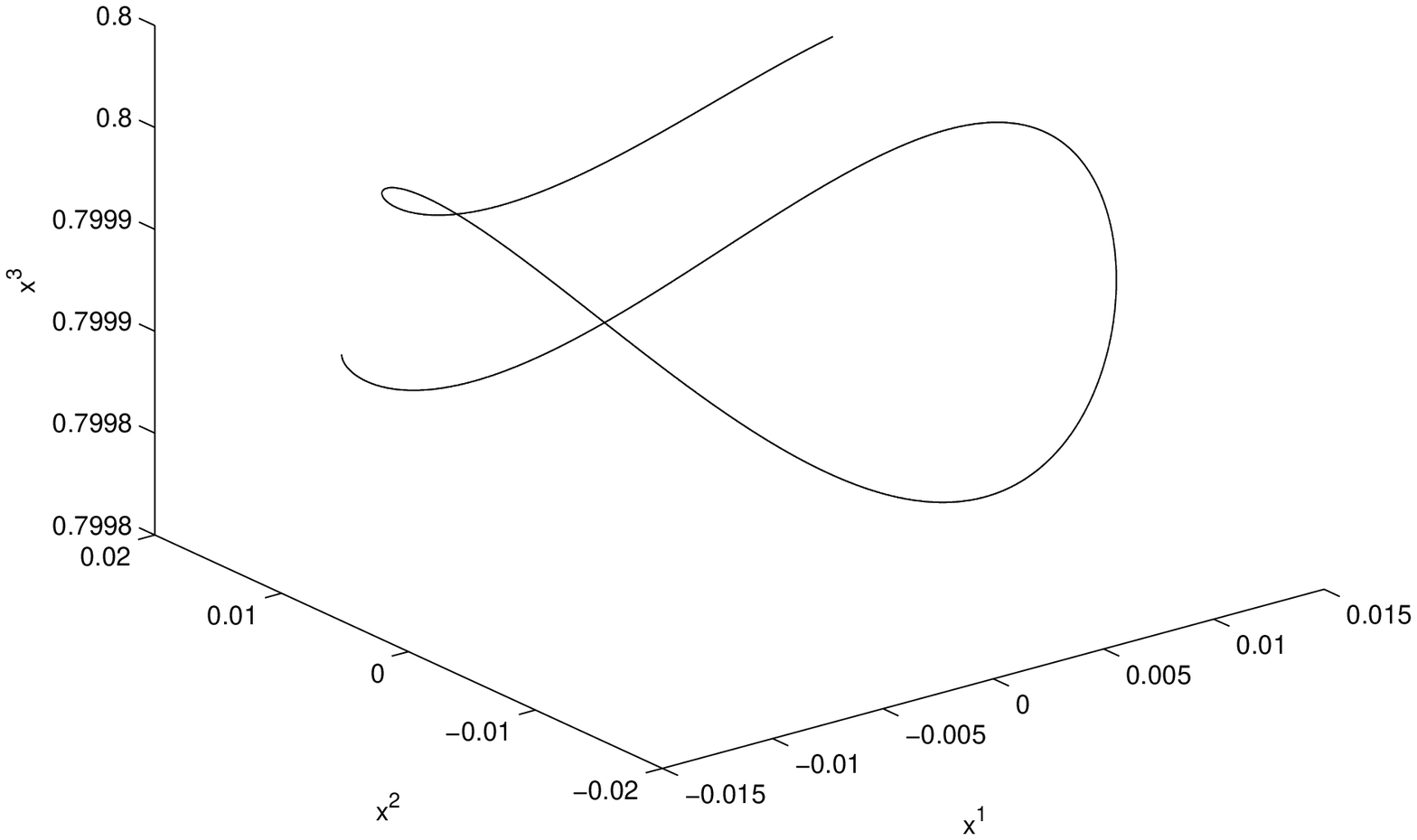} &
\epsfxsize=6cm \epsfysize=5cm \epsffile{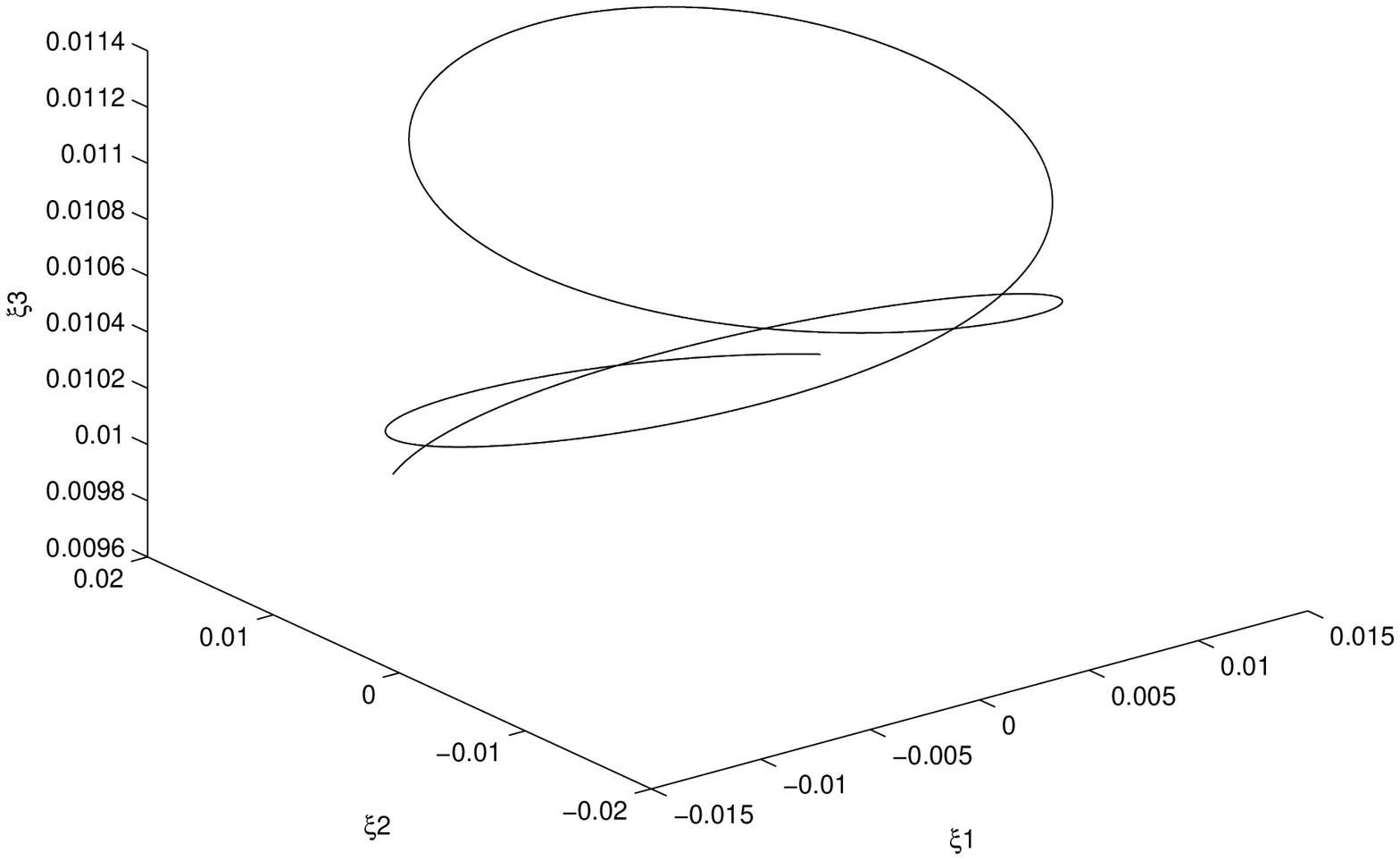}\\
Fig. 5 & Fig. 6
\end{tabular}
\end{center}

{\bf Conclusion}. In this paper are presented differential systems
on Leibniz algebroids. The analysis of these systems realizes by
studying the behavior of its solutions in the neighborhood of
equilibria points. This analysis will be enterprised in the future
papers.

\begin{center}
{\bf References}
\end{center}

{\bf [1]. D. Fish}, {\it Metriplectic systems}. Dissertation -
Portland State University, 2005.

{\bf [2]. J. Grabowski, P. Urbanski}, {\it Lie algebroids and Poisson - Nijenhuis structures}.
 Rep. Math. Phys., {\bf 40}(1997), 195 - 208.

{\bf [3]. J. Grabowski, P. Urbanski}, {\it Tangent and cotangent lifts and graded Lie algebras associated with Lie
algebroids}. Ann. Global Anal. Geom.,{\bf 15}(1997), no.5, 447 -486.

{\bf [4] J.L. Loday}, {\it Une version non-commutative des algebr$\grave e $s de Lie}. L'Einseignement Math$\acute e $ matique, {\bf 39} ( 1993), 269-293.

{\bf [5]. J. - P. Ortega, V. Planas -Bielsa}, {\it Dynamics on Leibniz manifolds}. Preprint
arXiv:math. DS/0309263,2003.

{\bf [6]. A. Weinstein} { In : Mechanics day ( Waterloo, ON,1992)} Fields Institute Communications 7, American Mathematical Society, 1996, 207. \\

Author's address:\\

West University of Timi\c soara\\
\hspace*{0.7cm} Department of Mathematics \\
\hspace*{0.7cm} 4, Bd. V. P{\^a}rvan, 300223, Timi\c soara, Romania\\
\hspace*{0.7cm} E-mail : ivan@hilbert.math.uvt.ro\\

West University of Timi\c soara\\
\hspace*{0.7cm} Department of Mathematics\\
\hspace*{0.7cm}4, Bd. V. P{\^a}rvan, 300223, Timi\c soara, Romania\\
\hspace*{0.7cm} E-mail : opris@math.uvt.ro\\
\end{document}